\documentclass{article}

\usepackage{amsfonts,amssymb,amsmath}

\newcommand{\reff}[1]{(\ref{#1})}
\usepackage{amsfonts,amssymb,amsmath}
\usepackage[T2A]{fontenc}
\usepackage[cp1251]{inputenc}
\usepackage{latexsym}
\usepackage{graphicx}

\newcommand{\ola}[1]{\overline{\alpha}}

\def\bbP{\mathbb P}
\def\bbR{\mathbb R}
\def\bbZ{\mathbb Z}
\def\bfm{\boldsymbol}
\def\cl{\centerline}

\begin{document}
\cl{\Large{\bf A non-linear model of limit order book dynamics}}
\vskip 5 truemm

\cl{{\bf N. Vvedenskaya}\footnote{Institute for Information
Transmission Problems, Russian Academy of Sciences, 19 Bol'shoi
Karetny Per, GSP-4, Moscow 127994, RUSSIA; {\it E-mail}
ndv$@$iitp.ru },
 {\bf Y. Suhov}
 \footnote{Institute for Information Transmission Problems,
Russian Academy of Sciences, 19 Bol'shoi Karetny Per, GSP-4, Moscow
127994, RUSSIA; Instito de Matematica e Estatistica, Universidade de
Sao Paulo, Caixa Postal 66281- CEP 05389-970 Sao Paulo, BRASIL;
Statistical Laboratoty, DPMMS, University of Cambridge, Wilberforce
Road, Cambridge CB3 0WB, UK; {\it E-mail} yms@statslab.cam.ac.uk} ,
{\bf V. Belitsky}\footnote{Instito de Matematica e Estatistica,
Universidade de Sao Paulo, Caixa Postal 66281- CEP 05389-970 Sao
Paulo, BRASIL; {\it E-mail} belitsky@ime.usp.br}}

 \vskip 5 truemm

\begin{abstract}This paper focuses on some simple models of limit order book
dynamics which simulate market trading mechanisms. We start with a
discrete time/space Markov process and then perform a re-scaling
procedure leading to a deterministic dynamical system controlled by
non-linear ODEs. This allows us to introduce approximants for the
equilibrium distribution of the process represented by fixed points
of deterministic dynamics. \vskip 5 truemm
\end{abstract}

\cl{\bf 1. Introduction. The underlying Markov process}
\vskip 5 truemm
\def\Mu{\rm M}
\def\Nu{\rm N}
\def\diy{\displaystyle}
\def\dist{{\rm{dist}}}

In what follows, LOB stands for the limit order book, a trading
mechanism adopted in many modern financial markets. For a detailed
description of some common LOB models and their applications, see
\cite{R} and references therein. (Although our models differ in a number
of aspects.) One of challenging problems is to determine factors
attracting or repelling various market participants.

This paper explores a new approach to the analysis of LOB dynamics
where the parameters of the original random (Markov) process are
re-scaled, and a limiting dynamical system emerges, with a
deterministic behavior described by a system of non-linear
(ordinary) differential equations. A similar approach
is commonly used in the literature on stochastic
communication networks; see, e.g., 
the paper \cite{VDK} and its sequels (in particular, \cite{VS}).
In the current paper we consider a simplified
model where a number of technically involved issues are absent. We
also omit proofs, referring the reader to 
forthcoming publications beginning with \cite{VSB}.

The rationale for the models below is as follows. We consider
a single-commodity market where prices may be at one of $N$ distinct
levels (say, $c_1<c_2<\ldots <c_N$, although the exact meaning of
these values is of no importance). The market is operating in
discrete times $0$, $\delta$,
$2\delta$, $\ldots$. At a given time $t\delta$, $t=0,1,\ldots $,
there are $b_i(t)$ traders
prepared to buy a unit of the commodity at price $c_i$ and $s_i(t)$ traders
prepared to sell it at this price, which leads to vectors

%
%
%
\begin{equation}\label{1}
 \bfm{b}(t)=\big(b_1(t),\ldots ,b_N(t)\big),\;
\bfm{s}(t)=\big(s_1(t),\ldots ,s_N(t)\big)\in\bbZ_+^N.
\end{equation}
Here $\bbZ_+=\{0,1,\ldots\}$ stands for a non-negative integer
half-lattice and $\bbZ_+^N$ for the non-negative integer $N$-dimensional
lattice orthant. The pair $(\bfm{b}(t), \bfm{s}(t))$ represent a
state of a Markov process $\bfm{U}(t)$ that will be the subject of
our analysis.

If $b_i(t)\geq s_i(t)>0$ then each of the sellers gets a trade with
probability $p_{\rm T}\in (0,1)$ and leaves the market, together
with his buyer companion. Therefore, both values $b_i(t)$ and
$s_i(t)$ decrease by a random number $n=0,1,\ldots, s_i(t)$ with the
binomial probability. A seller among $s_i-n$ who did not get the
trade either (i) quits the market with probability $p_{\rm Q}\in
(0,1)$ or (ii) moves to the price level $c_{i-1}$ with probability
$p_{\rm M} \in (0,1)$ or (iv) remains at the same level with
probability $1-p_{\rm Q} -p_{\rm M}$. (One can think that for this
seller a random experiment is performed, with three outcomes.)
Similarly, a buyer among $b_i(t)-n$ who did not get the trade quits
the market with probability $p_{\rm Q}\in (0,1)$ or moves to the
price level $c_{i+1}$ with probability $p_{\rm M}\in (0,1)$ or
remains at the same level with probability $1-p_{\rm Q} -p_{\rm M}$.
(Assuming that $p_{\rm Q} +p_{\rm M}< 1$.)

Symmetrically, if $s_i(t)\geq b_i(t)>0$ then each of the buyers gets
a trade with probability $p_{\rm T}\in (0,1)$ and leaves the market,
together with his seller companion. The remaining traders at the
price level $c_i$ proceed as above.


In addition, at every time point $t\delta$ a random number of new
buyers arrive and position themselves at the price level $c_1$; it
is distributed according to a Poisson law with mean $\Lambda_{\rm
b}>0$. Similarly, at every time $t\delta$ a Poisson random number of
sellers arrive and take a position at price level $c_N$; the mean
value of this variable is $\Lambda_{\rm s}>0$.

All described events occur at each level independently. This generates
the aforementioned Markov process $\big\{\bfm{U}(t)\big\}$ with trajectories
$\big\{(\bfm{b}(t),\bfm{s}(t))\big\}$, $t\in\bbZ_+$.
\vskip 5 truemm

{\bf Theorem 1.} {\sl $\forall$ values of parameters $\Lambda_{{\rm
b}/ {\rm s}}$, $p_{{\rm Q}/{\rm M}}$ and $p_{\rm T}$, the process
$\{\bfm{U}(t)\}$ is irreducible, aperiodic and
positive recurrent. Therefore, it has a unique set of equilibrium
probabilities
$\pi=\Big(\pi\big(\bfm{b},\bfm{s}\big):\;\bfm{b},\bfm{s}\in\bbZ^N\Big)$,
and $\forall$ initial state $\bfm{U}(0)$ (deterministic or random),
the distribution of the random state $\bfm{U}(t)$ at time $t$ converges
weakly to $\pi$ as $t\to\infty$:
\begin{equation}\label{2}
\lim_{t\to\infty}\bbP\Big(\bfm{U}(t)=(\bfm{b},\bfm{s})\Big)
=\pi\big(\bfm{b},\bfm{s}\big).
\end{equation}
}
\vskip 5 truemm
 \cl{\bf 2. Scaling limit}
 \vskip 5 truemm

The explicit form of the equilibrium distribution $\pi$ of process
$\big\{\bfm{U}(t)\big\}$ (and even
probabilities of transitions $(\bfm{b},\bfm{s})\mapsto (\bfm{b}',\bfm{s}')$)
are too cumbersome to work with. This fact makes it desirable to develop
efficient methods of approximation. In this paper we focus on one such
method based on scaling the parameters of the process (including states
and time-steps).

The re-scaling procedure is as follows: we set
\begin{equation}\label{3}
p_{\rm T}=\frac{\gamma}{L},\;\;\;p_{\rm Q}
=\frac{\alpha_{\rm Q}}{L},\;\;\;p_{\rm M}
=\frac{\alpha_{\rm M}}{L},\;\;\;
\Lambda_{\rm b}=\frac{\lambda_{\rm b}}{L},
\;\;\;\Lambda_{\rm s}=\frac{\lambda_{\rm s}}{L},
\end{equation}
 where $\gamma >0$, $\alpha_{\rm Q}>0$, $\alpha_{\rm M}>0$,
$\lambda_{\rm b}>0$
and $\lambda_{\rm s}>0$ are fixed and
$L\to\infty$. In addition, we re-scale the states and the time:
pictorially,
$$x_i\sim \frac{b_i}{L},\;\;y_i\sim \frac{s_i}{L},\;\;
\tau\sim\frac{t\delta}{L}.$$
Formally, denoting the Markov process generated for a given $L$
by $\bfm{U}^{(L)}$, we consider the continuous-time process
\begin{equation}\label{4}
\bfm{V}^{(L)}(\tau )={\diy\frac{1}{L}}\bfm{U}^{(L)} \big(\lceil\tau
L/\delta\rceil\big),\;\;\tau\geq 0, \end{equation}
 where $\lceil
a\rceil$ stands for the integer part of $a>0$.

Set: $\bbR_+=(0,\infty )$ (a positive half-line), then $\bbR_+^N$ is
a positive orthant in $N$ dimensions. Suppose we are given a pair of
vectors $(\bfm{x}(0),\bfm{y}(0))\in\bbR_+^N\times\bbR_+^N$ where
$\bfm{x}(0)=(x_1(0),\ldots ,x_N(0))$, $\bfm{y}(0)=(y_1(0),\ldots ,
y_N(0))$. Consider the following system of first-order ODEs for
functions $x_i=x_i(\tau )$ and $y_i=y_i(\tau )$ where $\tau >0$ and
$1\leq i\leq N$ :
\begin{equation}\label{5}
\begin{array}{l}
{\dot x}_1=\lambda_{\rm
b}-\Big(\alpha_{\rm Q}+ \alpha_{\rm M}\Big)x_1
-\gamma\min\;\big[x_1,y_1\big],\\
{\dot x}_i=\alpha_{\rm M}x_{i-1}
-\Big(\alpha_{\rm Q}+
\alpha_{\rm M}\Big)x_i
-\gamma\min\;\big[x_i,y_i\big],\;1<i\leq N,\\
{\dot y}_i=\alpha_{\rm M}y_{i+1}
\,-\Big(\alpha_{\rm Q}+
\alpha_{\rm M}\Big)y_i
\;-\gamma\min\;\big[x_i,y_i\big],\;1\leq i< N,\\
{\dot y}_N=\lambda_{\rm s}-\Big(\alpha_{\rm Q}+ \alpha_{\rm
M}\Big)y_N -\gamma\min\;\big[x_N,y_N\big],\end{array}
\end{equation}
 with
the initial date $x_i(0)$, $y_i(0)$, $1\leq i\leq N$. The fixed
point $\big(\bfm{x}^*,\bfm{y}^*\big)$ of system \reff{5} has
$\bfm{x}^*=(x^*_1,\ldots ,x^*_N)$ and $\bfm{y}^*=(y^*_1,\ldots
,y^*_N)$ where $x^*_i$ and $y^*_i$ give a solution to
\begin{equation}\label{6}
\begin{array}{l} \lambda_{\rm b}=\Big(\alpha_{\rm Q}+
\alpha_{\rm M}\Big)x^*_1
+\gamma\min\;\big[x^*_1,y^*_1\big],\\
\alpha_{\rm M}x^*_{i-1}=\Big(\alpha_{\rm Q}+\alpha_{\rm M}\Big)x^*_i
+\gamma\min\;\big[x^*_i,y^*_i\big],\;1<i\leq N,\\
\alpha_{\rm M}y^*_{i+1}
=\,\Big(\alpha_{\rm Q}+\alpha_{\rm M}\Big)y^*_i
+\;\gamma\min\;\big[x^*_i,y^*_i\big],\;1\leq i< N,\\
\lambda_{\rm s}=\Big(\alpha_{\rm Q}+ \alpha_{\rm M}
\Big)y^*_N +\gamma\min\;\big[x^*_N,y^*_N\big].\end{array}
\end{equation}
 Both systems \reff{5} and \reff{6} are non-linear. However, the
non-linearity `disappears' at a local level which greatly simplifies
the analysis of these systems.

In Theorems 2 and 3 below, we use the distance generated by the Euclidean
norm in $\bbR^N\times\bbR^N$.
\vskip 5 truemm

{\bf Theorem 2.} {\sl $\forall$ initial date
$(\bfm{x}(0),\bfm{y}(0))\in\bbR_+^N\times\bbR_+^N$ there exists a
unique solution $\Big\{(\bfm{x}(\tau ),\bfm{y}(\tau )),\;\tau
>0\Big\}$ to system} (5). {\sl For this solution, $(\bfm{x}(\tau
),\bfm{y}(\tau )) \in\bbR_+^N \times\bbR_+^N$ $\forall$ $\tau >0$.
As $\tau\to\infty$, the solution approaches a fixed point, which
yields a unique solution to system} (6):
\begin{equation}\label{7}
\dist\Big[\big(\bfm{x}(\tau ),\bfm{y}(\tau )\big),
\big(\bfm{x}^*,\bfm{y}^*\big)\Big]\to 0.
 \end{equation}

\vskip 5 truemm

{\bf Theorem 3.} {\sl Suppose that the re-scaled initial states
${\diy\frac{1}{L}}\bfm{U}(0)$ tend to vector
$(\bfm{x}(0),\bfm{y}(0))\in\bbR_+^N\times\bbR_+^N$ in probability:
$\forall$ $\epsilon >0$,
\begin{equation}\label{8}
\lim_{L\to\infty}\bbP\left(\dist\left[{\diy\frac{1}{L}}\bfm{U}(0),
(\bfm{x}(0),\bfm{y}(0))\right]\,\geq\epsilon\right) =0.
\end{equation}
 Then,
$\forall$ $T>0$, the process $\Big\{\bfm{V}^{(L)}(\tau ),\; \tau\in
[0,T]\Big\}$ converges in probability to the solution
$\Big\{(\bfm{x}(\tau ),\bfm{y}(\tau )),\;0\leq\tau\leq T\Big\}$. That
is, $\forall$ $\epsilon >0$,}
\begin{equation}\label{9}
\lim_{L\to\infty}\bbP\left(\sup\left\{\dist\left[\bfm{V}^{(L)}(\tau
), (\bfm{x}(\tau),\bfm{y}(\tau))\right] \,:\;0\leq\tau\leq T\right\}
\geq\epsilon\right)=0.
\end{equation}

{\sl In particular, if $\bfm{x}(0)=\bfm{x}^*$ and $\bfm{y}(0)=\bfm{y}^*$
then}
\begin{equation}\label{10}
\lim_{L\to\infty}\bbP\left(\sup\left\{\dist\left[\bfm{V}^{(L)}(\tau
), \big(\bfm{x}^*,\bfm{y}^*\big)\right] \,:\;0\leq\tau\leq T\right\}
\geq\epsilon\right)=0.
\end{equation}
 {\sl Moreover, if process
$\big\{\bfm{U}(t),\;t\in\bbZ_+\big\}$ is in equilibrium then Eqn}
(10) {\sl holds true.} \vskip 5 truemm

\cl{\bf 3. Fixed points in the scaling limit. Concluding remarks}
\vskip 5 truemm

The approximation developed in Theorem 3 calls for an analysis
of solutions to (6). As follows from the middle equations in (6),
\vskip 5 truemm

{\bf Lemma 4.} {\sl The fixed-point entries satisfy
\begin{equation}\label{11}
x^*_1>\ldots
>x^*_N\;\hbox{ and }\;y^*_1<\ldots <y^*_N.
\end{equation}
 Consequently,
the parameter space $\bbR^5_+$ formed by $\gamma$, $\alpha_{{\rm
Q}/{\rm M}}$, and $\lambda_{{\rm b}/{\rm s}}$ is
partitioned into open domains where one of the following generic
patterns persists:}

(i) $x^*_N>y^*_N$, (ii) $x^*_1<y^*_1$, {\sl and}
(iii) $x^*_i>y^*_i$
{\sl for $i=1,\ldots,\ell$ and $x^*_i<y^*_i$ for $i=\ell +1,\ldots
,N$ where $1<\ell<N$. In each of these domains system} \reff{6} {\sl
is linear.} \vskip 5 truemm

Lemma 4 allows us to develop simple algorithms for calculating the
fixed point $\big(\bfm{x}^*,\bfm{y}^*\big)$ and analyze the
character of convergence in \reff{7}.

A particular algorithm for calculating
$\big(\bfm{x}^*,\bfm{y}^*\big)$ is based on the following recursion.
Set $x^{(0)}_i=0$ and let $y^{(0)}_i$ be the solution to the third and
the forth equations \reff{6} with $y^{(0)}_N=\lambda_s/(\alpha_{\rm
Q}+ \alpha_{\rm M})$. Next, let
$\big(\bfm{x}^{(k)},\bfm{y}^{(k)}\big)$, $k=1,2,\dots$ be the
solution to the system
\begin{equation*}\label{12}
\begin{array}{l}
\lambda_{\rm b}=\Big(\alpha_{\rm Q}+ \alpha_{\rm M}\Big)x^{(k)}_1
+\gamma\min\;\big[x^{(k)}_1,y^{(k-1)}_1\big],\\
\alpha_{\rm M}x^{(k)}_{i-1} =\Big(\alpha_{\rm Q}+
\alpha_{\rm M}\Big)x^{(k)}_i
+\gamma\min\;\big[x^{(k)}_i,y^{(k-1)}_i\big],\;1<i\leq N,\\
\alpha_{\rm M}y^{(k)}_{i+1} =\,\Big(\alpha_{\rm Q}+
\alpha_{\rm M}\Big)y^{(k)}_i
+\;\gamma\min\;\big[x^{(k)}_i,y^{(k)}_i\big],\;1\leq i< N,\\

\lambda_{\rm s}=\Big(\alpha_{\rm Q}+ \alpha_{\rm M}\Big)y^{(k)}_N
+\gamma\min\;\big[x^{(k)}_N,y^{(k)}_N\big].\end{array}
\end{equation*}

{\bf Lemma 5.} {\sl The inequalities $x^{(k)}_i\ >\ x^{(k-1)}_i$,\
$y^{(k)}_i\ <\
y^{(k-1)}_i$\ hold true $\forall \ i,k\geq 1$\ and,  values
$x^{(k)}_i$ are uniformly
bounded. Therefore, $\exists$ $\lim\limits_{k\to \infty}x^{(k)}_i$,
and $\lim\limits_{k\to\infty}y^{(k)}_i$ and these limits satisfy
the system} \reff{6}.
\vskip 5 truemm

We conclude with the
following remarks. \vskip 5 truemm

1. The current set-up admits straightforward generalisations to the
case where parameters $\gamma$ and $\alpha_{{\rm Q}/{\rm M}}$
depend on $i$, $0<i<N$ and on the trader type
(b/s). A more complicated case emerges
if parameters $\lambda_{{\rm b}/{\rm s}}$ become state-dependent.

2. There are several forms of convergence for which the assertion in
Theorem 3 holds true. The dynamical system \reff{5} itself gives
rise to a limiting process with interesting properties.

3. Another valid approximation for process $\bfm{U}(t)$ is a
diffusion approximation working on a different scale from that in
\reff{3}. \vskip 5 truemm

These topics are subject to forthcoming research. See \cite{VSB} and
subsequent publications.


\vskip 5 truemm

\begin{thebibliography}{20}

\bibitem{R} I. Rosu. A Dynamic model of the limiting order book.
//faculty.chicagoboot.edu/ioanid.rosu/research/limit.pdf.


\bibitem{VDK} N.D. Vvedenskaya, R.L Dobrushin, F.I. Karpelevich. A queueing
system with selection of the shortest of two queues: an asymptotical
approach. {\it Problems of Information Transmission}, {\bf 32},
1996, 15--27.

\bibitem{VS} N.D. Vvedenskaya, Y.M. Suhov, Multy access system with many
users: stability and metastability. {\it Problems of Information
Transmission}, {\bf 43}, 2007, No 3.

\bibitem{VSB} N. Vvedenskaya, Y. Suhov, V. Belitsky. Non-linear models of
the limit order book dynamics, I. {\it In preparation} (2011).
\end{thebibliography}
\end{document}